\documentclass{amsproc}

\usepackage{amsmath}

\usepackage{xcolor}

 \usepackage{amssymb}

\usepackage{url}
\usepackage{supertabular} 

\usepackage{mathabx}
\usepackage{multirow} 
\usepackage{graphics}

\usepackage{float}

\usepackage{graphicx}

\allowdisplaybreaks 

\newcommand{\ud}{\mathrm{d}}    

\renewcommand{\Im}{\operatorname{Im}}
\newcommand{\sn}{\mathrm{sn}}
\newcommand{\cn}{\mathrm{cn}}
\newcommand{\dn}{\mathrm{dn}}

\renewcommand{\phi}{\varphi}

\theoremstyle{definition}

\theoremstyle{remark}

\theoremstyle{remark}

\theoremstyle{conjecture}

\numberwithin{equation}{section}

\title
{Elliptic Functions}
\date{\today}
\author{Shaun Cooper}
\address{
Institute of Mathematical and Computational Sciences, Massey University,
Private Bag 102904, North Shore Mail Centre, Auckland, New Zealand
E-mail: s.cooper@massey.ac.nz
}

\subjclass[2020]{Primary---33E05; Secondary---11F11, 35K05. \\To appear in ``Srinivasa Ramanujan: His Life, Legacy, and Mathematical Influence'', Springer,~2025}

\usepackage{hyperref}

\hypersetup{
  colorlinks = true,
  urlcolor = blue,
  linkcolor = black,
  citecolor = blue
}

\begin{document}

\begin{abstract}
This note discusses elliptic functions in Ramanujan's work.
\end{abstract}

\maketitle

An elliptic function is a function of a complex variable that is meromorphic and doubly periodic.
In general, Ramanujan does not use the standard notation for elliptic functions, nor does he mention
or make use of the double periodicity. Ramanujan's results are either best stated in his own
notation or written explicitly in terms of series and products. For example, central to Ramanujan's work
 on elliptic functions is the identity~\cite[(17)]{ramanujan_arithmetical}
\begin{align}
\label{Ram_id_1}
\lefteqn{ \left( \frac{1}{4} \cot \frac{\theta}{2} + \sum_{n=1}^\infty
\frac{q^n}{1-q^n}\sin n\theta \right)^2} \\
      \nonumber \\
& =  \left(\frac{1}{4} \cot \frac{\theta}{2} \right)^2
    + \sum_{n=1}^\infty \frac{q^n}{(1-q^n)^2} \cos n \theta
    + \frac{1}{2} \sum_{n=1}^\infty \frac{nq^n}{1-q^n}(1- \cos n \theta) \nonumber
\end{align}
where $q=\exp(2\pi i\tau)$ and $ |\Im\theta|<2\pi \Im\tau.$
It plays a significant role in the derivation
of Ramanujan's differential equations for Eisenstein series~\cite[(30)]{ramanujan_arithmetical} which are fundamental results in the theory of modular forms.
In his proof of~\eqref{Ram_id_1}, Ramanujan expands the left hand side and directly computes the coefficients in the Fourier
expansion to obtain the right hand side. Hardy liked Ramanujan's proof so much that he included it in his books~\cite{hardy} and~\cite{hw}.
Later, van der Pol~\cite{vanderpol} showed that Ramanujan's identity is a natural consequence of Jacobi's triple product identity together with the heat equation
$$
\frac{\partial u}{\partial t} = \frac{\partial^2 u}{\partial x^2}
$$
that is satisfied by the theta function. For fixed $\tau$, the function in parentheses on the left hand side of~\eqref{Ram_id_1}, namely
$$
f(\theta) = \frac{1}{4} \cot \frac{\theta}{2} + \sum_{n=1}^\infty\frac{q^n}{1-q^n}\sin n\theta,
$$
has an analytic continuation to a meromorphic function of $\theta$ whose singularities are
simple poles at each point of the lattice $\Lambda = \left\{2\pi m + 2\pi n \tau | m,n \in \mathbb{Z}\right\}.$
This underlying structure is responsible for the modular transformation properties of $f(\theta)$.  
Clearly $f(\theta+2\pi) = f(\theta)$, and it can be shown that $f(\theta+2\pi\tau) = f(\theta)-\frac{i}{2}$. It follows that the derivative $f'(\theta)$ is an
elliptic function of~$\theta$ with poles of order 2 which is closely connected to the Weierstrass~$\wp$ function.
Ramanujan's identity~\eqref{Ram_id_1} has powerful generalizations due to Venkatachaliengar~\cite[(1.21), (3.17)]{kv} that are studied further in~\cite[Thms. 1.16, 1.42]{cooper}.

One of the deepest results of Jacobi for elliptic functions is the inversion theorem, which may be stated as follows. Suppose $0<q<1$ and $x_4=x_4(q)$ is defined by
\begin{equation}
\label{jit1}
x_4=\left(\frac{\displaystyle{\sum_{n=-\infty}^\infty q^{(n+\frac12)^2}}}{{\displaystyle{\sum_{n=-\infty}^\infty q^{n^2}}}}\right)^4.
\end{equation}
Then $x_4$ increases from $0$ to $1$ as $q$ increases from $0$ to $1$, hence
the inverse function $q=q(x_4)$ exists. An explicit formula for the inverse function is given by
\begin{equation}
\label{jit3}
q=\exp\left( -\pi \, \frac{{}_2F_1(\frac12,\frac12;1;1-x_4)}{{}_2F_1(\frac12,\frac12;1;x_4)}\right)
\end{equation}
where ${}_2F_1$ is the hypergeometric function. Moreover,
$$
{}_2F_1\left(\frac12,\frac12;1;x_4\right) = \left(\sum_{n=-\infty}^\infty q^{n^2}\right)^2.
$$
Ramanujan formulated his own version of these results and sketched a proof~\cite[Ch.~17, Entries 2--6]{notebooks}
that has been completed in~\cite[pp. 91--102]{Part3}.
Ramanujan~\cite[pp. 257--262]{notebooks} went further to discover three incredible analogues of Jacobi's result. To state one of the analogues,
suppose $0<q<1$ and let $x_3=x_3(q)$ be defined by
\begin{equation}
\label{jit4}
x_3=\left(\frac{\displaystyle{\sum_{m=-\infty}^\infty\sum_{n=-\infty}^\infty q^{(m+\frac13)^2+(m+\frac13)(n+\frac13)+(n+\frac13)^2}}}
{{\displaystyle{\sum_{m=-\infty}^\infty\sum_{n=-\infty}^\infty q^{m^2+mn+n^2}}}}\right)^3.
\end{equation}
Then $x_3$ increases from $0$ to $1$ as $q$ increases from $0$ to $1$, hence
the inverse function $q=q(x_3)$ exists. An explicit formula for the inverse function is given by
\begin{equation}
\label{jit6}
q=\exp\left( -\frac{2\pi}{\sqrt{3}} \, \frac{{}_2F_1(\frac13,\frac23;1;1-x_3)}{{}_2F_1(\frac13,\frac23;1;x_3)}\right),
\end{equation}
and moreover
$$
{}_2F_1\left(\frac13,\frac23;1;x_3\right) = \sum_{m=-\infty}^\infty\sum_{n=-\infty}^\infty q^{m^2+mn+n^2}.
$$
Jacobi's example in~\eqref{jit1} and \eqref{jit3} is nowadays understood as being connected with modular forms of level~4.
Ramanujan's example in~\eqref{jit4} and \eqref{jit6} corresponds to the level~3 theory, while his other two analogues
correspond to levels~$1$ and~$2$. All four theories were employed by Ramanujan to construct his famous rapidly converging series for $1/\pi$ in~\cite{ramanujan_pi}.
A pioneering attempt to study Ramanujan's results on elliptic functions
was carried out by Venkatachaliengar~\cite{kv}. Significant contributions to the level~3 theory were made
by the Borwein brothers~\cite{borwein}. All of Ramanujan's results in~\cite[pp. 257--262]{notebooks} were proved in the seminal work of
Berndt, Bhargava and Garvan~\cite{bbg}. A different analysis has been given in~\cite[Chapters 3 and 4]{cooper} and the results have been extended to levels 5--12 in Chapters~5--12,
respectively, of~\cite{cooper}.

The Jacobian elliptic functions are introduced by Ramanujan in Chapter 18 of his second notebook~\cite{notebooks} where they are defined
by their Fourier series, namely 
\begin{align*}
S(\theta) &= \sum_{n=0}^\infty \frac{\sin(\frac12(2n+1)\theta)}{\sinh(\frac12(2n+1)y)}, \\
C(\theta) &= \sum_{n=0}^\infty \frac{\cos(\frac12(2n+1)\theta)}{\cosh(\frac12(2n+1)y)}, \\
\intertext{and}
C_1(\theta) &= \frac12+\sum_{n=1}^\infty \frac{\cos n\theta}{\cosh ny}.
\end{align*}
These are scaled versions of the Jacobian elliptic functions sn, cn and dn, respectively, and the precise identifications
with the Jacobian elliptic functions are given in~\cite[p. 168]{Part3}, \cite[p. 127]{cooper} or~\cite[p. 124]{kv}.
Ramanujan gives the differentiation formulas 
$$
\frac{\ud S}{\ud\theta} = CC_1, \quad \frac{\ud C}{\ud\theta} = -SC_1\quad\text{and}\quad \frac{\ud C_1}{\ud\theta} = -SC
$$ 
and offers analogues of the formulas
$$
\cn^2u+\sn^2u=1\quad\text{and}\quad \dn^2u+k^2\sn^2u = 1
$$
in terms of the functions $C$, $S$ and $C_1$. Ramanujan's ${}_1\psi_1$ summation formula can be used to
express each of the functions $C$, $S$ and $C_1$ as an infinite product, \cite[p. 126]{cooper}.

It is believed that Ramanujan did
not use complex function theory or double periodicity to study elliptic functions, and there remains
an air of mystery about how he may have made his discoveries~\cite[p. 212]{hardy}.
While Ramanujan's methods will likely never be known,
Venkatachaliengar's work~\cite{kv} offers a plausible and elegant approach to elliptic functions using methods that could have been available to Ramanujan.


\end{document}